\documentclass{amsart}
\usepackage{amsmath}
\usepackage{amsfonts}
\usepackage{amsthm}

\newtheorem{theorem}{Theorem}
\newtheorem{definition}{Definition}
\newtheorem{corollary}{Corollary}
\DeclareMathOperator{\Rep}{Rep}

\title{On $*$-representations of a certain class of algebras related
  to a graph}
\author{Vasyl Ostrovskyi}
\address{Institute of Mathematics, Nat. Acad. Sci. of Ukraine,
3, Tereschenkivska st., 01601, Kyiv, Ukraine}
\email{vo@imath.kiev.ua}
\date{\today}
\thanks{This research was partially supported by the State Foundation
  for Fundamental Research of Ukraine, grant no. 01.07/071 and by the
DFG, grant no.\ 436UKR 113/71/0-1 }

\begin{document}

\begin{abstract}
We study families of self-adjoint operators with given spectra whose
sum is a scalar operator. Such families are $*$-representations of
certain algebras which can be described in terms of graphs and
positive functions on them. The main result is that in the cases where the
graph is one of the extended Dynkin graphs $\tilde D_4$, $\tilde E_6$,
$\tilde E_7$ or $\tilde E_8$, all irreducible $*$-representations of the
corresponding algebra  are finite-dimensional. To prove this fact, we
introduce the notion of invariant functional on a graph and give their
description. 
\end{abstract}
\maketitle

\section*{Introduction}
Let $H$ be a separable Hilbert space, and let $A_1$, \dots, $A_n$ be a
family of self-adjoint operators in $H$ with fixed finite spectra,
such that $A_1+\dots+A_n=\lambda I$ for some $\lambda \in \mathbb
R$. Such family of operators can be treated as a representation of
certain $*$-algebra with finite number of generators and polynomial
relations. The corresponding algebras were studied in a number of
recent papers (see \cite{vla_mel_sam05} and the bibliography therein). The
interest to such algebras is due to their relations with deformed
preprojective algebras (see, e.g., \cite{cra_hol98}), representations
of quivers (see \cite{kru_roi05}), Horn problem (see \cite{kru_pop_sam05}),
integral operators (see \cite{vas98}) etc.

Considering a family of self-adjoint operators with a pre-defined
spectra, for which $A_1+\dots+A_n=\lambda I$ we can assume that
$\lambda>0$, and
\[
\sigma(A_l)\subset M_l = \{0=\alpha_0^{(l)} < \alpha_1^{(l)} < \dots
  <\alpha_{k_l}^{(l)}\},\quad  l=1, \dots, n.
\]  
To study such properties of operators, it is convenient to consider
them as $*$-representations of a certain $*$-algebra.
Following \cite{vla_mel_sam05} (see also \cite{mel_sam_zav04} and
references therein) consider a simply-laced non-oriented
graph $\Gamma$ consisting of $n$ 
branches, such that $l$-th branch has $k_l+1$ vertices,  $l=1$,
\dots, $n$, and all branches are connected at a single root vertex.
Marking the vertices of $l$-th branch (excluding the root vertex)  by
positive numbers $(\alpha_j^{(l)})_{j=1} ^{k_l}$ increasing to the
root,  we get a function 
\[
\chi=(\alpha_1^{(1)},\dots,\alpha_{k_1}^{(1)};\dots;
\alpha_1^{(n)},\dots,\alpha_{k_n}^{(n)})  
\] 
on the graph $\Gamma$ defined in
all veritces except the root (below this function will be called a character
on  $\Gamma$). The root vertex
will be marked by the 
number $\lambda$ (notice that the term character in many papers is
used to denote a function $(\chi,\lambda)$ on the whole graph,
including the root vertex. For our needs the given notation is more
convenient).  

Given a graph $\Gamma$, a character $\chi$ on $\Gamma$ and a positive
number $\lambda$, on can construct the following $*$-algebra
\[
\mathcal A_{\Gamma,\chi,\lambda}=\mathbb C\Bigl\langle a_l=a_l^*,
l=1,\dots,n \mid  p_l(a_l)=0, l=1,\dots,n; \ \sum_{l=1}^n a_l =\lambda e
\Bigr \rangle, 
\] 
where $p_l(x)=x(x-\alpha_1^{(l)})\dots(x-\alpha_{k_l}^{(l)})$, $k=1$,
\dots, $n$.  Then the family $A_1$, \dots, $A_n$ is a
$*$-representation of $\mathcal A_{\Gamma,\chi,\lambda}$.

Properties of the algebra $\mathcal A_{\Gamma,\chi,\lambda}$, in
particular, the structure of its $*$-reprsenta\-ti\-ons, crucially
depend on the type of the graph $\Gamma$: they are quite different for
the cases where $\Gamma$ is a Dynkin graph, extended Dynkin graph or
none of them. 

If $\Gamma$ is a Dynkin graph then the corresponding algebra is
finite-dimensional and therefore has only finite number of irreducible
$*$-representations, and all of them are finite-dimensional. In this
case, sets of parameters 
$(\chi,\lambda)$, for which 
there exist representations, and the $*$-representations themselves are
described in \cite{kru_pop_sam05,sam_zav05}. 

If $\Gamma$ is an extended Dynkin graph, then the algebra is
infinite-dimensional of polynomial growth \cite{vla_mel_sam05}. In
this case, there 
exists a special character $\chi_\Gamma$ on the corresponding extended
Dynkin graph $\Gamma$.
%, such that $(TS)^m:(\chi_\Gamma,\lambda) \mapsto
%(\chi_\Gamma,  \lambda^{(m)})$; here $m=1$, $2$, $3$ and $5$ for the
%extended Dynkin graphs $\tilde D_4$, $\tilde E_6$, $\tilde E_7$ and
%$\tilde E_8$ respectively. 
For such special characters, sets
$\Sigma_{\Gamma, \chi_\Gamma}$ of those $\lambda$, for which there
exist representations of the corresponding algebra and the
corresponding $*$-representations were  were studied in 
\cite{ost_sam99,mel_rab_sam04,mel_sam_zav04} etc. In particular, it
follows that for the special characters, all irreducible
$*$-representations are finite-dimensional.
% For special 
%characters, this set is infinite and contains unique limit point
%$\omega_\Gamma$. 

In Section~3 we prove that all irreducible
$*$-representations of $\mathcal A_{\Gamma,\chi,\lambda}$ are
finite-dimensional for any $\chi$ and $\gamma$ provided that $\Gamma$
is an extended Dynkin 
graph. As a corollary we get that for $\lambda\ne \omega(\chi)$ the
rigidity index for irreducible representation is equal to 2. 

To study the properties of $*$-representations of $\mathcal
A_{\Gamma,\chi,\lambda}$  in the case of arbitrary 
characters and more general graphs we introduce the notion of
invariant functional on the set of characters. We show that such
invariant functional is unique if and only if the graph is an extended
Dynkin graph; for Dynkin graph there are no invariant functionals, and
for other graphs there are exactly two invariant functionals
(Section~2). These invariant functionals are described in terms of
solutions of 
the equation 
\[
n-s = \sum_{l=1}^n \frac{1}{1+(s-1) + \dots + (s-1)^{k_l}}, \quad s\ge 1.
\]
In Section~1 we show that this equation has unique solution $s=2$ if
and only if $\Gamma$ is an extended 
Dynkin graph. For ordinary Dynkin graph this equation has no solutions
$s\ge 1$, and for other graphs it has two solutions $1<s_1<2<s_2<n$.

The main tool to study $*$-representations of $\mathcal
A_{\Gamma,\chi,\lambda}$ are reflection (Coxeter)  functors
introduced in \cite{kru02} (for non-involutive case, see
\cite{gel_pon71}). Namely 
there exist two functors, $S\colon \Rep \mathcal
A_{\Gamma,\chi,\lambda} \to \Rep \mathcal A_{\Gamma,\chi',\lambda'}$
and $T\colon \Rep \mathcal
A_{\Gamma,\chi,\lambda} \to \Rep \mathcal
A_{\Gamma,\chi'',\lambda}$, where 
\begin{align*}
\chi'&=(\alpha_{k_l}^{(1)}-\alpha_{k_l-1}^{(1)},\dots,
\alpha_{k_l}^{(1)}-\alpha_{0}^{(1)}; \dots; \alpha_{k_l}^{(n)}-
\alpha_{k_l-1}^{(n)},\dots, \alpha_{k_l}^{(n)}-\alpha_{0}^{(n)}),
\\
\lambda' &= \alpha_{k_1}^{(1)} + \dots +\alpha_{k_n}^{(n)} -\lambda,
\\
\chi''&=  (\lambda-\alpha_{k_1}^{(1)},\dots, \lambda -\alpha_1^{(1)};
\dots ;\lambda-\alpha_{k_n}^{(n)},\dots, \lambda -\alpha_1^{(n)}).
\end{align*} 
The action of these functors on $*$-representations gives rise to the
action on pairs, $S\colon (\chi,\lambda) \mapsto (\chi',\lambda')$,
$T\colon (\chi,\lambda )\mapsto (\chi'',\lambda)$. This action is
extensively used below.

\section{Equation which distinguishes the type of a graph}
Let $n$ and $k_l$, $l=1$, \dots, $n$ be given natural numbers, and
let $\Gamma$ be the corresponding graph. 
In what follows, we will use solutions of the following equation
\begin{equation}\label{graph:eq_with_s}
n-s = \sum_{l=1}^n \frac{1}{1+(s-1) + \dots + (s-1)^{k_l}}, \quad s\ge 1.
\end{equation}

\begin{theorem}\label{th:numbers} 
The equation \eqref{graph:eq_with_s} have no solutions on $[1,\infty)$
if and only if 
the corresponding graph $\Gamma$ is one of the Dynkin graphs $A_d$,
$d\ge1$, $D_d$, $d \ge 4$, $E_6$, $E_7$, or $E_8$.
The equation \eqref{graph:eq_with_s} has a unique solution on
$[1,\infty)$ if and only 
if the corresponding graph $\Gamma$ is one of the extended Dynkin
graphs $\tilde D_4$, $\tilde E_6$, $\tilde E_7$, or $\tilde E_8$,
this solution is $s=2$. 
In all other cases (i.e., where $\Gamma$ is neither a Dynkin graph
nor an extended Dynkin graph), the equation \eqref{graph:eq_with_s}
has on $[1,\infty)$ two solutions $1<s_1<2<s_2<n$.
\end{theorem}

\begin{proof}
Consider auxiliary functions $\phi_k(x)=(1+x+\dots+x^k)^{-1}$, so that
\eqref{graph:eq_with_s} takes the form 
\[
f_\Gamma(s)=n-s - \sum_{l=1}^n \phi_{k_l}(s-1) =0.
\]
Calculations give the following formula
\[
\phi_k''(x)=(k+1)\, x^{k-1}\,\bigl (kx^{k-1} + 2(k-1)x^{k-2} +\dots
+ (k-1) 2 x +k\bigr)\,\phi_k^3(x),
\]
which implies that $f_\Gamma''(s) <0$, $s>1$.

For each of the extended Dynkin graphs, a direct calculation shows
that ${f_\Gamma(2)  =0}$ and $f'_\Gamma(2) =0$, which gives the
uniqueness of the solution for the case of extended Dynkin graphs.

Each of the Dynkin graphs $D_4$, $E_6$, $E_7$, $E_8$ is a subgraph of
the corresponding extended Dynkin graph, which gives inequalities for
the corresponding functions, $f_{D_4}(s)<f_{\tilde D_4}(s)$,
$f_{E_6}(s)<f_{\tilde E_6}(s)$, $f_{E_7}(s)<f_{\tilde E_7}(s)$,
$f_{E_8}(s)<f_{\tilde E_7}(s)$ for $s\ge1$. For the graphs $A_n$, $n \ge
1$, and $D_n$, $n>4$, one checks directly the inequality
$f_\Gamma(s)<0$, $s\ge 0$, where $\Gamma$ is one of these graphs.

Let $\Gamma$ be neither Dynkin graph, nor extended Dynkin graph. Then $\Gamma$
contains subgraph $\Gamma_0$ which is an extended Dynkin graph, and
therefore, the inequality $f_\Gamma(s)> f_{\Gamma_0}(s)$ holds for
$s\ge 1$, in particular, $f_\Gamma(2)>0$. To complete the proof notice
that $f_\Gamma(1)=-1$ and $f_\Gamma(n)<0$, and $f''(s)<0$ guarantees
that there are only two solutions on the interval $[1,\infty)$. 
\end{proof}

\section{Invariant functionals on graphs}
Let $\Gamma$ be a graph formed by $n$ branches connected in a single
root vertex, and let $k_l$ be the number of vertices (excluding root)
in $l$-th branch. Let $\chi$ be a character on the graph $\Gamma$,
\begin{gather}\label{chi}
\chi = (\alpha_1^{(1)}, \dots, \alpha_{k_1}^{(1)}; \dots;
\alpha_1^{(n)},
\dots,\alpha_{k_1}^{(n)}),
\\
0<\alpha_1^{(l)}<\dots
<\alpha_{k_l}^{(l)},\quad
l=1, \dots, n. \notag
\end{gather}
Let $\omega(\cdot)$ be a linear functional, which takes non-negative
values on characters. 

\begin{definition}
We say that $\omega(\cdot)$ is invariant with respect to the functor
$TS$, if   
\[
TS(\chi,\omega(\chi)) = (\tilde\chi,\omega(\tilde\chi))
\]
for any character $\chi$ on $\Gamma$.
\end{definition}

\begin{theorem}
If $\Gamma$ is a Dynkin graph (one of $D_d$, $d\ge4$, $E_6$, $E_7$ or
$E_8$), then there are no invariant functionals on the set of its
characters.

If  $\Gamma$ is an extended Dynkin Graph (one of $\tilde D_4$,
   $\tilde E_6$, $\tilde E_7$ or $\tilde E_8$), then there exists unique
   invariant functional:

--- for $\tilde D_4$, $\omega(\alpha;\beta;\gamma;\delta) =
    \frac12(\alpha+\beta +\gamma +\delta)$;

--- for $\tilde E_6$,
    $\omega(\alpha_1,\alpha_2;\beta_1,\beta_2;\gamma_1,\gamma_2) =
    \frac13(\alpha_1 +\alpha_2 +\beta_1 +\beta_2 +\gamma_1 +\gamma_2)$;

--- For $\tilde E_7$,
    $\omega(\alpha_1,\alpha_2,\alpha_3;\beta_1,\beta_2
    ,\beta_3;\gamma) = \frac14(\alpha_1 +\alpha_2 +\alpha_3 +\beta_1 +\beta_2
    +\beta_3 +2\gamma)$;

--- for $\tilde E_8$, $\omega(\alpha_1,\alpha_2,
    \alpha_3,\alpha_4,\alpha_5;\beta_1,\beta_2;\gamma) =  \frac16(\alpha_1
    +\alpha_2 +\alpha_3+\alpha_4 + \alpha_5 + 2 \beta_1 +2\beta_2
    +3\gamma)$. 

In the case where $\Gamma$ is neither Dynkin graph nor extended
   Dynkin graph, there exist two $TS$-invariant functionals. They are
   given by
\begin{equation}\label{omega_general}
\omega(\chi) = \sum_{l=1}^n \sum_{j=1}^{k_l} a_j^{(l)} \alpha_j^{(l)},
\quad a_j^{(l)} \ge 0, \quad j=1,\dots, k_l; \ l=1,\dots,n.
\end{equation}
with
\begin{equation}\label{ajl}
a_j^{(l)}=\frac{(s-1)^j} {1+(s-1)+\dots+(s-1)^{k_l}}, \quad
j=1,\dots,k_l; \quad l=1,\dots, n,
\end{equation}
where $s$ is a solution of \eqref{graph:eq_with_s}.
\end{theorem}

\begin{proof}
Let $\chi$ be given by \eqref{chi}, then linear functional
$\omega(\cdot)$ can be represented as \eqref{omega_general}. 
Then $(TS)(\chi,\omega(\chi)) =
(\tilde\chi,\tilde\lambda)$, where
\begin{align*}
\tilde\chi&= (\sum_{l=1}^n \alpha_{k_l}^{(l)} - \sum_{l=1}^n \sum
_{j=1}^{k_l} a_j^{(l)} \alpha_j^{(l)} -
\alpha_{k_1}^{(1)},
\\
&\qquad \sum_{l=1}^n \alpha_{k_l}^{(l)} - \sum_{l=1}^n \sum 
_{j=1}^{k_l} a_j^{(l)} \alpha_j^{(l)} - \alpha_{k_1}^{(1)}
+\alpha_1^{(1)},  \dots,
\\
&\qquad \sum_{l=1}^n \alpha_{k_l}^{(l)} - \sum_{l=1}^n \sum 
_{j=1}^{k_l} a_j^{(l)} \alpha_j^{(l)} - \alpha_{k_1}^{(1)}
+\alpha_{k_1-1}^{(1)};  \dots),
\\
\tilde\lambda&= \sum_{l=1}^n \alpha_{k_l}^{(l)} - \sum_{l=1}^n \sum
_{j=1}^{k_l} a_j^{(l)} \alpha_j^{(l)}.
\end{align*}
Then 
\begin{align*}
\omega(\tilde\chi)-\tilde\lambda&= \sum_{l=1}^n \bigl( \alpha_1^{(l)}
((1-s)a_1^{(l)} +a_2^{(l)}) 
+ \dots + \alpha_{k_l-1}^{(l)}((1-s) a_{k_l-1}^{(l)} + a_{k_l}^{(l)})\notag 
\\
&\qquad \quad + \alpha_{k_l}^l( s-s_l -1 + (1-s)a_{k_l}^l)\bigr) 
\end{align*}
where $s=\sum_{l=1}^n s_l $,
$s_l=\sum_{p=1}^{k_l} a_p^{(l)}$. Since the latter should be zero for
arbitrary choice of $\chi$, we get the following conditions
\begin{gather}
a_{k_l}^{(l)}= (s-1) a_{k_l-1}^{(l)},\quad \dots,\quad  a_2^{(l)} =
(s-1)a_1^{(l)}, \label{aux:akl}
\\
a_{k_l}^{(l)} = 1- \frac{s_l}{s-1},\quad l=1,\dots,n,\label{aux:akls}
\end{gather} 
For $s_l$, $l=1$, \dots, $n$ from \eqref{aux:akl} we have
\begin{equation}\label{aux:sl}
s_l = a_1^{(l)}+\dots + a_{k_l}^{(l)} =a_1^{(l)}(1+(s-1)+\dots+(s-1)^{k_l-1}) 
\end{equation}
Substitute in \eqref{aux:akls} $a_{k_l}^{(l)}=(s-1)^{k_l-1} a_1^{(l)}$, then
$(s-1)^{k_l}a_1^{(l)} =s -1-s_l$, or 
\begin{equation}\label{aux:al1}
a_1^{(l)} = \frac{s-1}{1+(s-1)+\dots+(s-1)^{k_l}}
\end{equation}
Compare this with \eqref{aux:sl} we get
\begin{equation}\notag
s_l=
1-\frac{1} {1+(s-1)+\dots+(s-1)^{k_l}}
\end{equation}
And taking into account that $s=\sum_{l=1}^n s_l $ we finally get a
conditions \eqref{graph:eq_with_s} for $s$. Solutions of this
equations are described by Theorem~\ref{th:numbers}. 
For such $s$, using \eqref{aux:akl}, \eqref{aux:al1}, we get
expression \eqref{ajl} for $a_j^{(l)}$.
\end{proof}

\section{Irreducible representations}
Let $\Gamma$ be an extended Dynkin graph, and let $\chi$ be a
character on it. The main result on representations of the
corresponding algebra $A_{\Gamma,\chi,\lambda}$ is the following.

\begin{theorem}
All irreducible families of operators corresponding to extended
   Dynkin diagrams are finite-dimensional.
\end{theorem}

\begin{proof}
Let $\pi$ be an irreducible representation of the algebra
$A_{\Gamma,\chi,\gamma}$, where $\Gamma$ is an extended Dynkin
graph. We consider two cases.

1. Let $\lambda=\omega(\chi)$. It is shown in
   \cite{mel05,vla_mel_sam05} that the 
   corresponding algebra is finite-dimensional over its center, and
   therefore, is a PI-algebra. This implies that the dimensions of all
   its irreducible representations are bounded.    

2. Let $\lambda<\omega(\chi)$. We proceed as follows. We apply  the
   $(ST)^n$ functors to  the representation of the algebra corresponding
   to the pair $(\chi,\lambda)$ to get representations of the  algebras
   corresponding to other pairs $(\chi_n,\lambda_n)$ and show that at
   some step either there cannot exist representation (in this case,
   there are no representations for $(\chi,\lambda)$), or the
   representation is an obvious extension of representation  of a
   subgraph (such subgraph is a Dynkin diagram, and the corresponding
   algebra is finite dimensional, therefore it has a finite number of
   representations all of which are finite-dimensional). Then the
   initial representation $\pi$ is obtained from some
   finitely dimensional representation of $A_{\Gamma,\chi_n,\lambda_n}$
   as a result of applying of the $(TS)^n$ functor, and therefore, is
   finite dimensional as well. 

Notice first that if some of the coefficients of $\chi$ is greater
than $\lambda$, then the corresponding projection is zero. Indeed, let
$P_j$, $j=1$,\dots,$m$ be orthogonal projections such that
$\sum_{j=1}^m \alpha_j P_j=\lambda I$, where $\alpha_j>0$,
$j=1$,\dots,$m$, and $\alpha_k>\lambda$ for some $k$. Then $\sum_{j\ne
  k} \alpha_jP_kP_jP_k =(\lambda-\alpha_k)P_k$ which is possible only
for $P_k=0$ since the left-hand side of the latter equality is
non-negative, but the right-hand side is non-positive.

In the case where some of the coefficients of $\chi$ are equal to
$\lambda$, the corresponding projection commutes with all other
projections and therefore, is either identity (in this case all other
projections are zero), or zero.

Thus, in both these cases, the representation is in fact a
representation of subalgebra in $A_{\Gamma,\chi,\lambda}$ corresponding
to some subgraph.

Now let all coefficients of $\chi_k$, $k\le n$ be positive, and some
coefficient of $\chi_{n+1} $ be negative or zero. Taking into account the way
the functors act on characters, we easily see that this means that the
corresponding coefficient of $\chi_n$ is grater or equal that $\lambda_n$.

To complete the proof, it is now sufficient to show that for any
$\lambda <\omega(\chi)$ there exists such $n$, that some  coefficient
of $\chi_n$ is negative or zero.

Let $\chi_\Gamma$ be the special character of the corresponding graph,
and let $\omega_\Gamma=\omega(\chi_\Gamma)$. If we norm $\chi$ such that
$\omega(\chi)=\omega_\Gamma$, then the character can be represented as
$\chi = \chi_\Gamma + \tilde \chi$, where $\tilde \chi$ is (not necessarily
positive) character, such that $\omega( \tilde \chi) =0$. Also
represent $\lambda= \omega_\Gamma-\lambda$. That is:

--- for $\tilde D_4$ $\chi=(1+a_1; 1+a_2; 1+a_3;1+a_4)$,
    $\lambda=2-\gamma$, $a_1+a_2+a_3+a_4=0$;

--- for $\tilde E_6$ $\chi=(1+a_1,2+a_2;1+b_1,2+b_2;1+c_1,2+c_2)$,
    $\lambda=3-\gamma$, $a_1+a_2+b_1+b_2+c_1+c_2=0$; 

--- for $\tilde E_7$
    $\chi=(1+a_1,2+a_2,3+a_3;1+b_1,2+b_2,3+b_3;2+c_1)$,
    $\lambda=4-\gamma$, $a_1+a_2+a_3+b_1+b_2+b_3+2c_1=0$ ;

--- for $\tilde E_8$ $\chi=(1+a_1, 2+a_2, 3+ a_3, 4+ a_4, 5+a_5;
    2+b_1,4+b_2; 3+c_1)$, $\lambda=6-\gamma$,
    $a_a+a_2+a_3+a_4+a_5+2b_1+2b_2 +3c_1=0$.

Using this notation, one can directly check that for any $k=1$, 2, \dots,

--- for $\tilde D_4$ $(ST)^{2k} (\chi,\lambda) =(\chi - 2k\gamma
    \chi_\Gamma,\lambda - 4k\gamma) $;

--- for $\tilde E_6$ $(ST)^{6k}(\chi,\lambda) = (\chi -3k \gamma
    \chi_\Gamma, \lambda-9k \gamma)$;  

--- for $\tilde E_7$ $(ST)^{12k} (\chi,\lambda) =(\chi - 4k\gamma
    \chi_\Gamma,\lambda - 16k\gamma) $;

--- for $\tilde E_8$ $(ST)^{30k} (\chi,\lambda) =(\chi - 6k\gamma
    \chi_\Gamma,\lambda - 36k\gamma)$.

Here $\chi_\Gamma$ is the special character for the corresponding
graph. The relations listed above complete the proof in the case
$\lambda <\omega(\chi)$. 

3. Let $\lambda >\omega(\chi)$. Apply  the $S$ functor to the pair
   $(\chi,\lambda))$, then we get a pair $(\chi',\lambda')$ with
   $\lambda'<\omega(\chi')$ and the arguments above apply.
\end{proof}

Recall that the rigidity index \cite{kat96,str_voe99} of a family of
operators $A_j$, 
$j=1$, \dots, $k$ in $n$-dimensional space is
\[
r=n^2(2-k) + \sum_{j=1}^k c(A_j)
\]
where $c(A_j)$ is the dimension of a centralizer of $A_j$. 

\begin{corollary}
The rigidity index is equal to 2 for all irreducible representations
of $A_{\Gamma,\chi,\gamma}$ where $\Gamma$ is extended Dynkin graph
and $\lambda\ne \omega(\chi)$.
\end{corollary}

\begin{proof}
Indeed, rigidity index is preserved by $T$ and $S$ functors. From the
proof above it follows that any irreducible representation can be
obtained from one-dimensional ones under the action of the Coxeter
functors (in fact, it was shown that any irreducible representation is
obtained from representation of some algebra related to ordinary
Dynkin graph, but this process can be iterated for the subalgebra
until we finish in one-dimensional space). Now the result follows from
the directly verified fact that for 
one-dimensional representations $r=2$.
\end{proof}

The author  expresses his gratitude to Professor
Yu.~S.~Samoilenko for kind attention and fruitful discussions of the
topics discussed in this paper. 

%\bibliography{all}

\begin{thebibliography}{VMS05}

\bibitem[CBH98]{cra_hol98}
W.~Crawley-Boevey and M.~P. Holland, \emph{Noncommutative deformations of
  {K}leinian singularities}, Duke Math. J. \textbf{92} (1998), no.~3, 605--635.

\bibitem[GP71]{gel_pon71}
I.~M. Gel'fand and V.~A. Ponomarev, \emph{Quadruples of subspaces of a
  finite-dimensional vector space}, Dokl. Akad. Nauk SSSR \textbf{197} (1971),
  no.~4, 762--765, (Russian).

\bibitem[Kat96]{kat96}
N.~Katz, \emph{Rigid local systems}, Princeton Univ. Press, 1996.

\bibitem[KPS05]{kru_pop_sam05}
S.~Kruglyak, S.~Popovych, and Yu. Samoilenko, \emph{Representations of algebras
  associated with {D}ynkin graphs and spectral problem}, Journ. Algebra Appl.
  (2005), to appear.

\bibitem[KR05]{kru_roi05}
S.~A. Krugljak and A.~V. Roiter, \emph{Locally scalar representations of graphs
  in the category of {H}ilbert spaces}, nct. Anal. Prilozh.F \textbf{39}
  (2005), no.~2.

\bibitem[Kru02]{kru02}
S.~A. Kruglyak, \emph{Coxeter functors for a certain class of $*$-quivers and
  $*$-algebras}, Methods Funct. Anal. Topol. \textbf{8} (2002), no.~4, 49--57.

\bibitem[Mel05]{mel05}
A.~Mellit, \emph{Certain examples of deformed preprojective algebras and
  geometry of their $*$-representations}, arXiv:math.RT/0502055 (2005).

\bibitem[MRS04]{mel_rab_sam04}
A.S. Mellit, V.I. Rabanovich, and Yu.S. Samoilenko, \emph{When a sum of partial
  reflections is a scalar operator}, Funct. Anal. Prilozh. \textbf{38} (2004),
  no.~2, 91--94.

\bibitem[MSZ04]{mel_sam_zav04}
A.~Mellit, Yu. Samoilenko, and M.~Zavodovsky, \emph{On $*$-representations of
  algebras of {T}emperley-{L}ieb type and algebras generated by linearly
  dependent generators with given spectra}, Proc. Inst. Math. NAS of Ukraine
  \textbf{50} (2004), 1139--1144.

\bibitem[OS99]{ost_sam99}
V.~L. Ostrovskyi and Yu.~S. Samoilenko, \emph{Introduction to the theory of
  representations of finitely presented $*$-algebras. {I}. {R}epresentations by
  bounded operators}, vol.~11, Rev. Math.\& Math. Phys., no.~1, Gordon \&
  Breach, London, 1999.

\bibitem[SV99]{str_voe99}
K.~Strambach and H.~V\"olklein, \emph{On linearly rigid tuples}, J. Reine
  Angew. Math. \textbf{510} (1999), 57--62.

\bibitem[SZ05]{sam_zav05}
S.~Samoilenko, Yu and M.~V. Zavodovsky, \emph{Spectral theorems for
  $*$-representations of the algebras {$\mathcal{P}_{\Gamma,\chi,com}$}
  associated with {D}ynkin graphs}, Methods Funct. Anal. Topol. \textbf{11}
  (2005), no.~1, 88--96.

\bibitem[Vas98]{vas98}
N.~Vasilevski, \emph{${C}^*$-algebras generated by orthogonal projections and
  their applications}, Integral Equations and Operator Theory \textbf{31}
  (1998), 113--132.

\bibitem[VMS05]{vla_mel_sam05}
M.~S. Vlasenko, A.~S. Mellit, and Yu.~S. Samoilenko, \emph{On algebras
  generated with linearly dependent generators that have given spectra}, Funct.
  Anal. Appl. \textbf{39} (2005), no.~3, to appear.

\end{thebibliography}
%\bibliographystyle{amsalpha}

\providecommand{\bysame}{\leavevmode\hbox to3em{\hrulefill}\thinspace}
\providecommand{\MR}{\relax\ifhmode\unskip\space\fi MR }
\providecommand{\MRhref}[2]{%
  \href{http://www.ams.org/mathscinet-getitem?mr=#1}{#2}
}
\providecommand{\href}[2]{#2}

\end{document}